\documentclass[twoside,leqno,twocolumn]{article}

\usepackage[letterpaper]{geometry}

\usepackage{ltexpprt}

\usepackage{pifont}

\usepackage{epsfig}
\usepackage{amssymb}
\usepackage{amsmath}
\usepackage{amsfonts}
\usepackage{bbding}
\usepackage{array}
\usepackage{caption,tabularx,booktabs}
\usepackage{arydshln}
\usepackage{paralist}
\usepackage{xargs}                      
\newcommand\tagthis{\addtocounter{equation}{1}\tag{\theequation}}

\newcommand{\m}[1]{{\bf{#1}}}

\newcommand{\C}[1]{{\cal {#1}}}

\usepackage{xcolor}

\usepackage{array}
\usepackage{caption,tabularx,booktabs}
\usepackage{arydshln}

\usepackage{xargs}                      
\usepackage{mathtools}

\newcommand{\g}[1]{\bm #1}

\newcommand{\R}{\mathbb{R}}

\newcommand{\sign}{\text{sgn}}

\usepackage{bm}

\newcommand{\norm}[1]{\left\|#1\right\|}
\usepackage{algorithm}
\usepackage{algorithmic}
\usepackage{subfigure}

\usepackage{enumitem}
\usepackage{multirow}
\def\argmin{\mathop{\rm argmin}}
\def\argmax{\mathop{\rm argmax}}

\usepackage{multirow}

\usepackage[font=small,skip=5pt]{caption}

\newtheorem{defn}{Definition}

\usepackage{bm}

\usepackage{algorithm}
\usepackage{algorithmic}
\usepackage{subfigure}

\usepackage{enumitem}
\usepackage{multirow}
\def\argmin{\mathop{\rm argmin}}
\def\argmax{\mathop{\rm argmax}}

\newcolumntype{C}[1]{>{\centering\let\newline\\\arraybackslash\hspace{0pt}}m{#1}}

\makeatletter
\newcounter{subthm} 
\let\savedc@thm\c@hyp

\newcommand{\normhyp}{%
  \let\c@hyp\savedc@hyp 
  \renewcommand\thehyp{\arabic{hyp}}%
} 
\makeatother

\makeatletter
\newcounter{subass} 
\let\savedc@ass\c@hyp

\makeatother

\usepackage[bottom]{footmisc}

\usepackage{fancyhdr} 
\begin{document}

\pagenumbering{arabic} 


\title{\Large StepDIRECT - A Derivative-Free Optimization Method for Stepwise Functions}
\author{Dzung T. Phan$^{1,}$\thanks{IBM Research, T. J. Watson Research
Center, New York, USA, Email: \{phandu@us.,lamnguyen.mltd@\}ibm.com} \and Hongsheng Liu$^{1,}$\thanks{Work done while an intern at IBM Research. University of North Carolina, Chapel Hill, NC, USA, E-mail: hsliu@live.unc.edu} \and
Lam M. Nguyen\footnotemark[1]}

\date{}

\maketitle

\def\thefootnote{1}\footnotetext{Equal contribution}\def\thefootnote{\arabic{footnote}}

\pagenumbering{arabic}







\begin{abstract} \small\baselineskip=9pt 
In this paper, we propose the \texttt{StepDIRECT} algorithm for derivative-free optimization (DFO), in which the black-box objective function has a stepwise landscape. Our framework is based on the well-known \texttt{DIRECT} algorithm. By incorporating the local variability to explore the flatness, we provide a new criterion to select the potentially optimal hyper-rectangles. In addition, we introduce a stochastic local search algorithm performing on potentially optimal hyper-rectangles to improve the solution quality and convergence speed. Global convergence of the \texttt{StepDIRECT} algorithm is provided. Numerical experiments on optimization for random forest models and hyper-parameter tuning are presented to support the efficacy of our algorithm. The proposed \texttt{StepDIRECT} algorithm shows competitive performance results compared with other state-of-the-art baseline DFO methods including the original \texttt{DIRECT} algorithm. 
\end{abstract}

\section{Introduction}\label{sec: intro}
We introduce an optimization algorithm for solving a class of \textit{structured} black-box deterministic problems, which often arise in  data mining and machine learning 
\begin{equation}\label{eq:original_problem_1}
 \min_{\m{x}\in {\Omega}}  f(\m{x}), 
\end{equation}
where ${\Omega}$ is a bounded hyper-rectangle, i.e., ${\Omega}=\{\m{x} \in \R^{p} : \m{l} \le \m{x} \le \m{u}\}$ for some given bounds $\m{l},\m{u} \in \R^{p}$. We assume that $f : \R^{p} \to \R$ is a \textit{stepwise} function, whose closed-form formula is unavailable or costly to get and store. 
A stepwise function is a piecewise constant function over a finite number of disjoint subsets. 
We will point out two motivating important applications in the next section that fit into this framework: hyper-parameter tuning (HPT) for classification and optimizing tree ensemble regression. 

Derivative-free optimization (DFO) has a long history and can be traced back to the deterministic direct-search (DDS) method proposed in \cite{hooke1961direct}. DFO algorithms can be classified into two categories: \textit{local} and \textit{global} search methods. Local algorithms focus on techniques that can seek a local minimizer. Direct local algorithms find search directions by evaluating the function value directly; for example, Nelder–Mead algorithm \cite{Nelder1965ASM} and the generalized
pattern-search method \cite{torczon1997convergence}. Model-based algorithms construct and optimize a local surrogate model for the objective function to determine the new sample point. For instance, the radial basis function is utilized in RBFOpt \cite{costa2018rbfopt} as the surrogate model, while polynomial models are used in \cite{Powell1994,Conn97recentprogress}. Due to the local flatness and discontinuous structure of the stepwise function in \eqref{eq:original_problem_1}, a local search algorithm might easily get stuck at a bad local minimizer. We note that \eqref{eq:original_problem_1} can have an infinite number of local minimizers in some applications and gradients are either zero or undefined. Global search algorithms aim at finding a global solution. Methods based on Lipschitzian-based partitioning techniques for underestimating the objective include the dividing  rectangles algorithm (DIRECT) \cite{jones1993lipschitzian} and branch-and-bound search \cite{Pinter97}. Stochastic search algorithms such as particle swarm optimization (PSO) \cite{Kennedy95} and the differential evolution (DE) \cite{price2006differential} are also considered. Global model-based approaches optimize a surrogate model, usually constructed for the entire search space, e.g., response surface methods \cite{Jones2001b}. These methods often require a large number of samples in order to obtain a high-fidelity surrogate. Recently, Bayesian optimization using Gaussian process \cite{frazier2018tutorial,mockus2012bayesian} are widely applied in black-box optimization, especially hyper-parameter tuning for machine learning models.


Due to the excessive number of local minimizers and the computational cost of a function evaluation in many applications such as the hyper-parameter tuning problem, we focus on a global optimization algorithm that can reduce function values quickly in a moderate number of iterations. In this work, we propose a global spatial-partitioning algorithm for solving the problem \eqref{eq:original_problem_1}, based on the idea of the well-known \texttt{DIRECT} algorithm \cite{jones1993lipschitzian}. The key ingredient in our algorithm is the selection of hyper-rectangles to do partition and sample new points.  We do not attempt to compute approximate gradients or build a surrogate of the objective function as in many existing methods. The name \texttt{DIRECT} comes from the shortening of the
phrase ``DIviding RECTangles”, which describes the way the algorithm partitions the feasible domain by a number of hyper-rectangles in order to move towards the optimum. An appealing feature of \texttt{DIRECT} is that it is insensitive to discontinuities and does not rely on gradient estimates. These characteristics are nicely suitable for solving the stepwise function~\eqref{eq:original_problem_1}.

There are two main steps in a \texttt{DIRECT}-type algorithm: 1) selecting a sub-region within $\Omega$, the so-called potentially optimal hyper-rectangle, in order to get a new sample point over the sub-region, and 2) splitting the potentially optimal hyper-rectangle.  In the literature, a number of variants of \texttt{DIRECT} algorithms have been proposed to improve the performance of \texttt{DIRECT}, most of them are devoted to the first step \cite{Finkel2006,Liuzzi2010,Jones2020}. There are a limited number of papers working on the second step, for example \cite{Grbic2013}. Another line of research direction for speeding up \texttt{DIRECT} is to incorporate a local search strategy in the framework \cite{Liuzzi2016}. Most of these modifications are for a general objective function, but to the best of our knowledge, very few of them have been successful to exploit the problem structure or prior knowledge on the objective function \cite{Grbic2013,Liuzzi2010}. For example, in \cite{Grbic2013}, the authors present an efficient modification of \texttt{DIRECT} to optimize a symmetric function by including an ordered set in the hyper-rectangle dividing step. In \cite{Liuzzi2010}, the authors assume that the optimal function value is known in the hyper-rectangle selection step. In this paper, we propose a new \texttt{DIRECT}-type algorithm to utilize the stepwise landscape of the objective function, and make contributions in both two steps. Compared with the original \texttt{DIRECT} and its variants, the proposed \texttt{StepDIRECT} differs from them in the following aspects: 
%
\vspace{-4mm}
\begin{itemize}
\item We provide a new criterion for selecting potentially optimal hyper-rectangles using the local variability and the best function value. 
\vspace{-2mm}
\item We propose a new stochastic local search algorithm, specially designed for stepwise functions, to explore the solution space more efficiently. As a result, \texttt{StepDIRECT} is a stochastic sampling global optimization algorithm, where incorporating stochasticity can help it to escape from a local minimum.
\vspace{-6mm}
\item When prior information on the relative importance for decision variables is available, we split the potentially optimal hyper-rectangles along the dimension with higher relative variable importance.
\end{itemize}
\vspace{-2mm}
We can prove that the algorithm converges asymptotically to the global optimum under mild conditions. Experimental results reveal that our algorithm performs better than or on par with several state-of-the-art methods in most of the settings studied.

\vspace{-1mm}
\section{Motivating Examples}\label{sec:motivating_examples}

To fully see the importance of the proposed algorithm, we now show two concrete examples. Before giving a detailed description of these examples, we present a formal definition of a stepwise function as follows.

\begin{defn}
A function $f : \R^p \rightarrow \R$ is called a stepwise function over $\Omega \subset \R^p$ if there is a partition $\{\Omega_i\}_{i=1}^D$ of $\Omega$ and real numbers $\{c_i\}_{i=1}^D$ such that 
$\Omega_i\cap \Omega_j=\emptyset 
$
for $1\leq i < j \leq D$ and $
f(\m{x}) = \sum_{i=1}^D c_i \mathbf{I}{(\m{x}\in \Omega_i)}$.  
Here,  $\m{I}(\cdot)$ is the indicator function.
\end{defn}

\subsection{Tree-based Ensemble Regression Optimization}
A common approach for the decision-making problem based on data-driven tools is to build a pipeline from historical data, to predictive model, to decisions. A two-stage solution is often used, where the prediction and the optimization are carried out in a separate manner \cite{DemirovicCPAIOR2019}. Firstly, a machine learning model is trained to learn the underlying relationship between the controllable variables and the output. Secondly, the trained model is embedded in the downstream optimization to produce a decision. We assume that the regression model estimated from the data is a good representation of the complex relationship. In this paper, we consider the problem of optimizing a tree ensemble model such as random forests and boosting trees, where the predictive model has been learned from historical data.

The tree ensemble model combines predictions from multiple decision trees.  A decision tree uses a tree-like structure to predict the outcome for an input feature vector $\m{x} \in \R^p$. The $t$-th regression tree in the ensemble model has the following form
\begin{equation}\label{eq:DecisionTree}
f_t(\m{x})=\sum_{i=1}^{M_t}c_{t,i} \cdot \mathbf{I}{(\m{x}\in \Omega_{t,i})}
\end{equation}
where $\Omega_{t,1},\ldots, \Omega_{t,M_t}$ represent a partition of feature space. We can see that $f_t(\m{x})$ is a stepwise function. 
The tree ensemble regression function  outputs predictions by taking the weighted sum of multiple decision trees as 
\begin{equation}\label{eq:TreeEnsemble}
f(\m{x})=\sum_{t=1}^{T} w_t f_t(\m{x}),
\end{equation}
where $w_t$ is the weight for the decision tree $f_t(\m{x})$. We assume that parameters $c_{t,i}, \Omega_{t,i}$ and $w_t$ have been learned from data, and we use \texttt{StepDIRECT} to optimize $f(\m{x})$.  As we can see that $f_t(\m{x})$ is constant over each sub-region $R_{t,i}$; hence $f(\m{x})$ is a stepwise function. Formally, we have the following theorem.  

\begin{theorem}
Regression functions for decision tree ensemble methods using a constant value for prediction at each leaf node including random forest and gradient boosting are stepwise functions.
\end{theorem}

We note that for this type of regression functions, we might get additional information about the objective function such as \textit{variable importance} for input features of random forests \cite{breiman2001random}.

\subsection{Hyper-parameter Tuning  for Classification}
In machine learning models, we often need to choose a number of hyper-parameters to get a high prediction accuracy for unseen samples \cite{Hutter2019AutomatedML}. For example, in training an $\ell_1$-regularized logistic regression classifier, we tune the sparsity parameter \cite{Lee2006EfficientLR}. 
A common practice is to split the entire dataset into three subsets: training data, validation data, and test data \cite{hastie2009elements}. For a given set of hyper-parameters, we train the model on the training data, then evaluate the model performance on the validation data. 
Some widely-used performance metrics include accuracy, precision, recall, $F_1$-score, and AUC (Area Under the Curve). The goal is to tune hyper-parameters for the model to maximize the performance on the validation data. The task can be formulated as a black-box optimization problem.  

\begin{figure}[h!]
    \centering
    \includegraphics[width=0.80\linewidth,height=0.20\textheight]{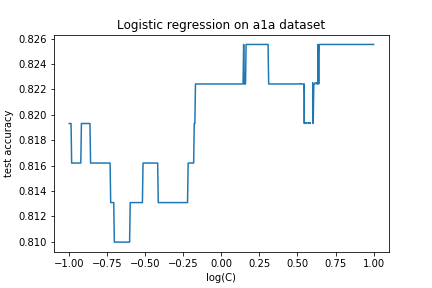}
    \caption{The landscape of the objective function of the HPT problem}
    \label{fig:htp_prob}
\end{figure}

First, we start with a binary classifier. Suppose that we are given the $M$ training samples $\{(\m{u}_1,\m{v}_1),\ldots,(\m{u}_M,\m{v}_M)\}$ and $N$ validation samples $\{(\m{u}_1,\m{v}_1),$ $\ldots, (\m{u}_N,\m{v}_N)\}$ where $\m{u}_i \in \R^{d}$ and $\m{v}_i \in \{\pm 1\}$. For a fixed set of model parameters $\g{\lambda} \in \R^{p}$, a classification model $h(\cdot;\g{\lambda}) : \R^{d} \to \{\pm 1\}$ has been learned based on the training data. When the performance measure is accuracy, the HPT problem is to determine $\g{\lambda}$ that maximizes the test accuracy
\[
F_{acc}(\g{\lambda}) = \frac{1}{N}\sum_{i=1}^N\m{I}(h(\m{u}_i;\g{\lambda}) = \m{v}_i),
\]
where $(\m{u}_i,\m{v}_i)$ is from the validation data. We can see that $F_{acc}(\g{\lambda}) \in \{0,\frac{1}{N},\ldots,\frac{N-1}{N},1\}$ for any $\g{\lambda}$; hence we expect that the landscape of this function should have stepwise behavior. As an example, we plot in Figure \ref{fig:htp_prob} the landscape for the HPT logistic regression by tuning the $\ell_1$ parameter (denoted by $C$ in the $x$-axis) for  training the \texttt{a1a} dataset \cite{CC01a}.  

The target function for other metrics can also take only a finite number of values. The observation still holds true for a multi-label classifier. 
\section{\texttt{StepDIRECT} for Stepwise Functions}\label{sec:direct}




The \texttt{StepDIRECT} algorithm begins the optimization by transforming the domain of the problem \eqref{eq:original_problem_1} linearly into the unit
hyper-cube. Therefore, we assume for the rest of the paper that
\begin{equation}
{\Omega} = \{\m{x}\in \R^p: 0\leq x_i \leq 1\}.
\end{equation} 
In each iteration, \texttt{StepDIRECT} consists of three main steps. First, we identify a set of \textit{potentially optimal} hyper-rectangles based on a new criterion. We expect that the hyper-rectangles have a high chance to contain a global optimal solution. The second step is to perform a local search over the potentially optimal hyper-rectangles. Thirdly, we divide the selected hyper-rectangles into smaller hyper-rectangles. 

We improve over the general \texttt{DIRECT}-type approaches by: 1) using a different heuristic to select which hyper-rectangle to split, which takes into account the local function variability; 2) using a local search (randomized directional search) to choose a high-quality point as a representative for each hyper-rectangle. Both two proposed strategies make use of the stepwise structure of the objective function. 


At the $k$-th iteration, let $\C{P}_k$ be the set of all generated hyper-rectangles $\C{H}_i$ in the partition of $\Omega$, where   $\C{H}_i = \{\m{x} \in \R^{p}:\m{l}_i \le \m{x} \le \m{u}_i\}$ for some bounds $\m{l}_i$ and $\m{u}_i$.
We let $f_{\C{H}_i}$ denote the best function value of $f$ over $\C{H}_i$ (by evaluating at the sampling points in the sub-region $\C{H}_i$). 
$f_{min}$ and $\m{x}_{min}$ are the best function value and the best feasible solution over $\Omega$, respectively. We use $m$ to count the number of function evaluations and $m^{max}$ is the maximal number of function evaluations. We present our main algorithm \texttt{StepDIRECT} in Algorithm~\ref{alg:stepdirect}.     

\begin{algorithm}[htb!]
	\caption{\texttt{StepDIRECT}}\label{alg:stepdirect}
	\begin{algorithmic}
	\STATE Define $\m{c}_1=(0.5, \ldots, 0.5)\in \R^p$ and set $\m{x}_{min}=\m{c}_1, f_{min}=f(\m{c}_1)$, $k=m=1$
	\STATE Run the Initialization Step to get $\C{P}_1, f_{\C{H}_i} \mbox{ for } \C{H}_i \in \C{P}_1, f_{min}$ and $\m{x}_{min}$  (see Sect. \ref{sec:init}) 
	\WHILE{$m\leq m^{max}$}
	\STATE Identify the set $\C{S}$ of all potentially optimal hyper-rectangles in $\C{P}_k$ (see Sect. \ref{sec:potentialRect}) 
	\FOR{$j\in \C{S}$}
	\STATE a) Perform a local search over $\C{H}_j$ to get $f_{\C{H}_j}$ (see Sect. \ref{sec:localsearch}) 
	\STATE b) Identify the sides of the rectangle $\C{H}_j$, divide $\C{H}_j$ into smaller hyper-rectangles along these sides, and update $\C{P}_k$  (see Sect. \ref{sec:dividing}) 
	\STATE c) Evaluate $f$ at centers of new hyper-rectangles, and update $f_{\C{H}_i}$ for every $\C{H}_i \in \C{P}_k$.  
	\ENDFOR
	 \STATE Update $m, f_{min}$ and $x_{min}$ from Steps a and c.
	 \STATE Update $k = k+1$ 
	\ENDWHILE
	\end{algorithmic}
\end{algorithm}


\subsection{Initialization Step}
\label{sec:init}

We follow the first step of the original \texttt{DIRECT} for initialization \cite{finkel2003direct}. In this step, Algorithm \ref{alg:stepdirect} starts with finding $f_1 = f(\m{c}_1)$ and divides the hyper-cube ${\Omega}$ by evaluating the function values at $2p$ points $\m{c}_1\pm\delta \m{e}_i, i\in  \{1,\ldots,p\}$, where $\delta=\frac{1}{3}$ and  $\m{e}_i$ is the $i$-th unit vector. The idea of \texttt{DIRECT} is to select a hyper-rectangle with a small function value in a large search space;
hence let us define
\[
s_i = \min\{f(\m{c}_1+\delta \m{e}_i), f(\m{c}_1 - \delta \m{e}_i)\}, \forall i\in \{1,\ldots,p\}
\]
and the dimension with the smallest $s_i$ is partitioned into thirds. By doing so, $\m{c}_1\pm\delta \m{e}_i$ are the center of the newly generated hyper-rectangles. We initialize the sets $\C{P}_1,\C{I}_1, \C{C}_1$, values $f_{\C{H}_i}$ for every $\C{H}_i \in \C{P}_1$, and update $f_{min} = \min\{f_{\C{H}_i} : \C{H}_i \in \C{P}_1\}$ and corresponding $\m{x}_{min}$. 




\subsection{Potentially Optimal Hyper-rectangles}
\label{sec:potentialRect}
In this subsection, we propose a new criterion for \texttt{StepDIRECT} to select the next potentially optimal hyper-rectangles, which should be divided in this iteration. In the original \texttt{DIRECT} algorithm, every hyper-rectangle $\C{H}_i$ is represented by a pair $(f_{\C{H}_i},d_i)$, where $f_{\C{H}_i}$ is the function value estimated at the centre $\m{c}_i$ of $\C{H}_i$ and $d_i$ is the distance from the center of hyper-rectangle $\C{H}_i$
to its vertices. The criterion to select hyper-rectangles, i.e., \textit{potentially optimal hyper-rectangles}, for further divided is based on a score computed from $(f_{\C{H}_i},d_i)$. A pure local strategy would select the hyper-rectangle with the smallest value for $f_{\C{H}_i}$, while a pure global search strategy would choose one of the hyper-rectangles with the biggest value $d_i$. The main idea of the \texttt{DIRECT} algorithm is to balance between the local and global search, which can be achieved by using a score weighting the two search strategies: $f_{\C{H}_i} - K d_i$ for some $K > 0$. The original definition for potentially optimal hyper-rectangles $\C{H}_j$  for \texttt{DIRECT} \cite{jones1993lipschitzian} is based on two conditions:
\begin{align} 
f(\m{c}_j) - K d_j \; &\leq \; f(\m{c}_i) - K d_i, \;\; \forall \C{H}_i\in \C{P}_k  \label{eq:poc_010} \\
f(\m{c}_j) - K d_j \; &\leq \; f_{min} - \epsilon, \label{eq:poc_020}
\end{align}
for some $\epsilon > 0$. 
%

\begin{figure}[h!]
    \centering
    \includegraphics[width=0.8\linewidth,height=0.22\textheight]{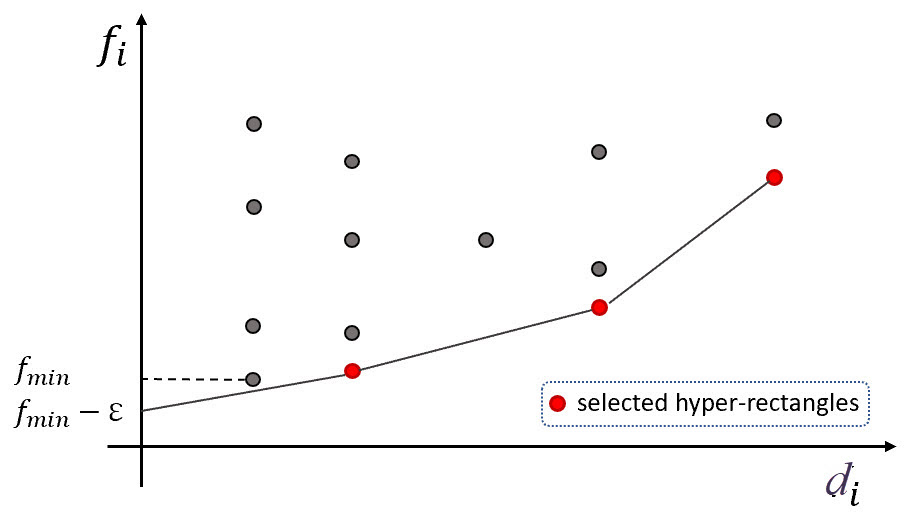}
    \caption{Identifying potentially optimal hyper-rectangles for \texttt{DIRECT}}
    \label{fig:pohr}
\end{figure}

An example for identifying these potentially optimal hyper-rectangles by the original \texttt{DIRECT} is given in Figure \ref{fig:pohr}. Each dot in a two dimensional space represents a hyper-rectangle, three red dots with the smallest value for $f(\m{c}_j) - K d_j$ for each $K$ and a significant improvement (i.e., $f(\m{c}_j) - Kd_j > f_{min} - \epsilon$) are considered as potentially optimal.

The proposed \texttt{StepDIRECT} searches locally and globally by dividing all hyper-rectangles that meet the criteria in Definition \ref{def:new_potential_optimal}. A hyper-rectangle $\C{H}_i$ is now represented by a triple $(f_{\C{H}_i},d_i,\sigma_i)$, in which we introduce new notations for $f_{\C{H}_i}$ and $\sigma_i$. We use a higher quality value $f_{\C{H}_i}$ returned by a local search as a representative for each hyper-rectangle, and a flatness measurement $\sigma_i$ in the proposed criterion.


\begin{defn} (potentially optimal hyper-rectangle) \label{def:new_potential_optimal}
Let $\epsilon>0$ be a small positive constant and $f_{min}$ be the current best function value over $\Omega$. Suppose $f_{\C{H}_i}$ is the best function value over the hyper-rectangle $\C{H}_i$. A hyper-rectangle $\C{H}_j$ is said to be potentially optimal if there exists $K>0$ such that
\begin{align} 
f_{\C{H}_j} - Kd_j \sigma_j \; &\leq \; f_{\C{H}_i} - K d_i \sigma_i, \quad \forall \C{H}_i\in \C{P}_k \label{eq:poc_1} \\
f_{\C{H}_j} - K d_j \sigma_j \; &\leq \; f_{min} - \epsilon |f_{min}-f_{median}|, \label{eq:poc_2}
\end{align}
where 
$\sigma_j > 0$ quantifies the local variability of $f$ on the hyper-rectangle $\C{H}_j$ (defined in \eqref{eq:vari}) and $f_{median}$ is the median of all  function values in history.
\end{defn}

In Eq. \eqref{eq:poc_1}, we introduce a new notation $\sigma_j$ to measure the local landscape of the stepwise function. Furthermore, we replace $f(\m{c}_j)$ as in the original \texttt{DIRECT} by $f_{\C{H}_j}$, which is computed with the help of local search. The value $f_{\C{H}_j}$ for $\C{H}_j$ generated during Steps a and c is better estimating the global solution of $f$ over $\C{H}_j$ because $f_{\C{H}_j} \le f(\m{c}_j)$. In Eq. \eqref{eq:poc_2}, in order to eliminate the sensitivity to scaling issue, $|f_{min}|$ is replaced by the difference  $|f_{min}-f_{median}|$ as suggested in \cite{Finkel2006}.

Now we discuss how to compute $\sigma_j$. At the $k$-th iteration of the \texttt{StepDIRECT} algorithm, for hyper-rectangle $\C{H}_j \in \C{P}_k$ with center $\m{c}_j$ and diameter $2d_j$, we define its set of \textit{neighborhood} hyper-rectangles as 
\[
\C{N}_j = \{\C{H}_i\in \mathcal{P}_k: \norm{\m{c}_i-\m{c}_j} \leq \lambda d_j\}
\]
for some $\lambda > 0$. 
Then, $\C{N}_j$ can be further divided into two disjoint subsets $\C{N}_j^D$ and $\C{N}_j^E$ such that
\begin{align*}
    \C{N}_j^D = \{\C{H}_i \in \C{N}_j: f_{\C{H}_i}\neq f_{\C{H}_j}\} \ \text{and} \\
    \C{N}_j^E = \{\C{H}_i\in \C{N}_j: f_{\C{H}_i}= f_{\C{H}_j}\}.
\end{align*}
The local variability  estimator for hyper-rectangle $\C{H}_j$ is defined as 
\begin{equation}\label{eq:vari}
\sigma_j = \max\Big\{\dfrac{|\C{N}_j^D| }{|\C{N}_j|},\epsilon_\sigma\Big\} \in [\epsilon_\sigma, 1],
\end{equation}
where $0<\epsilon_\sigma<1$ is a small positive number to prevent $\sigma_j=0$ when $\C{N}_j^D = \emptyset$. 
In our experiments, we set $\lambda=2$ and $\epsilon_\sigma=10^{-8}$ as default.
The meaning of $\sigma_j$ can be interpreted as follows. For each hyper-rectangle $\C{H}_j$, a large value for $\sigma_j$ indicates that it is likely to have more different function values in the neighbourhood of the center $\m{c}_j$, which requires more local exploration. We propose to include the hyper-rectangle to the set of potentially optimal hyper-rectangles.


By Definition \ref{def:new_potential_optimal}, we can efficiently find potentially optimal hyper-rectangles based on the following lemma.
\begin{lemma}\label{lem_0001}
Let $\epsilon>0$ be the positive constant used in Definition \ref{def:new_potential_optimal} and $f_{min}$ be the current best function value. Let $\mathcal{I}$ be the set of indices of all existing hyper-rectangles. For each $j\in \mathcal{I}$, let $
\mathcal{I}_l = \{i\in \mathcal{I}: d_i \sigma_i < d_j \sigma_j\}, 
\mathcal{I}_b = \{i\in \mathcal{I}: d_i \sigma_i > d_j \sigma_j\},
\mathcal{I}_e = \{i\in \mathcal{I}: d_i \sigma_i = d_j \sigma_j\}
$
and $g_{i} = \frac{f_{\C{H}_i} - f_{\C{H}_j}}{d_i \sigma_i - d_j \sigma_j}$ for all $i\neq j$.
The hyper-rectangle $\C{H}_j$ is potentially optimal if and only if three following conditions hold:
\begin{itemize}
    \item[(a)] 
$f_{\C{H}_j} \leq f_{\C{H}_i}$ for every $i \in \mathcal{I}_e;$
\item[(b)] 
$
\max_{i\in \mathcal{I}_l}g_{i}\leq  \min_{i\in \mathcal{I}_b}g_{i};
$
\item[(c)] If $f_{median} > f_{min}$ then 
\begin{align*}
    \epsilon \leq \frac{f_{min} - f_{\C{H}_j}}{|f_{min}-f_{median}|}  + \frac{d_j \sigma_j \min_{i\in \mathcal{I}_b}g_{i}}{|f_{min}-f_{median}|}; \tagthis \label{eq:lemma1_3}
\end{align*}
otherwise,  
\begin{equation}\label{eq:lemma1_4}
f_{\C{H}_j} \leq d_j \sigma_j \min_{i\in \mathcal{I}_b}g_{i} + f_{min}.
\end{equation}
\end{itemize}

\end{lemma}




We note that all proofs are delegated
to the supplementary material. Definition \ref{def:new_potential_optimal} differs from the potentially optimal hyper-rectangle definition proposed in \cite{jones1993lipschitzian} (i.e., \eqref{eq:poc_010},\eqref{eq:poc_020}) in the following three aspects:
\begin{itemize}
    \item \eqref{eq:poc_1} and \eqref{eq:poc_2} include the local variability of $f$ on the hyper-rectangle. For the stepwise function, this quantity helps balance the local exploration and global search.
\item \eqref{eq:poc_2} uses $|f_{min}-f_{median}|$ to remove sensitivity to the linear and additive scaling  and improve clustering of sample points near optimal solutions.
\item $f_{\C{H}_i}$ can be different from $f(\m{c}_i)$, i.e., $f_{\C{H}_i} \le f(\m{c}_i)$.
\end{itemize}

\subsection{Dividing Potentially Optimal Hyper-rectangles}
\label{sec:dividing}
Once a hyper-rectangle has been identified as potentially optimal, \texttt{StepDIRECT} divides this hyper-rectangle into smaller hyper-rectangles. We will take into account the side length  and the variable importance measure of each dimensions. 

In tree ensemble regression optimization, we can obtain the variable importance which indicate the relative importance of different features \cite{breiman2001random}. In general, the variable importance relates to a significant function value change if we move along this direction, and also the number of discontinuous points along the coordinate. As a result, we tend to make more splits along the direction with higher variable importance. Formally, let $\m{w} = (w_1, \ldots, w_p)\in \R_{+}^p$ be the normalized variable importance with $||\m{w}||_1=1$ and the length of the hyper-rectangle be $\m{l}=(l_1, \ldots, l_p)\in \R_{+}^p$, then we define 
\begin{equation}
\m{v} = (v_1, \ldots, v_p) = (w_1l_1,\ldots, w_pl_p)
\end{equation}
as the relative importance of all coordinates for the selected hyper-rectangle. We choose the coordinate with the highest value $v_i$ to divide the hyper-rectangle into thirds. 

If no variable importance measure is provided, then we take $\m{w} = (1/p, \ldots,$ $1/p)\in \R_{+}^p$. The dividing procedure is the same as the original \texttt{DIRECT} by splitting the dimensions with the largest side length.

\subsection{Local Search for Stepwise Function}
\label{sec:localsearch}
In this subsection, we introduce a randomized local search algorithm designed for minimizing stepwise function $f(\m{x})$ over a bounded box \[\C{B} = \{\m{x} \in \R^{p} : \m{a} \le \m{x} \le \m{b}\}.\] 
It has been known that a DIRECT-type algorithm has a good ability to locate promising regions of the feasible space, and a good local search procedure can help to quickly converge to the global optimum \cite{Liuzzi10,Liuzzi2016}. Since the function values over the neighborhood of a point are almost surely constant for a stepwise function, existing local searches based on a local approximation of the objective function might fail to move to a better feasible solution. Hence, we propose a novel local search method, shown in Algorithm \ref{alg:LS}. Different from the classical trust-region methods, Algorithm \ref{alg:LS} will increase the stepsize when no better solution is discovered in the current iteration, i.e., $\tau > 1$. This change is motivated by the stepwise landscape of the objective function. 


\begin{algorithm}[htb!]
	\caption{Local\_Search ($\C{B})$}\label{alg:LS}
	\begin{algorithmic}
	\STATE \textbf{input}  $\C{B}$
	\STATE Given starting point $\m{x}\in \C{B}$, $\tau > 1, 0< \delta_{min} < \delta_{max}, \delta\in [\delta_{min}, \delta_{max}]$, $t^{max} > 0$. 
	\STATE Initialize $t=0$, $f^{\C{B}}_{min}=f(\m{x})$, and  $\m{x}^{\C{B}}_{min} = \m{x}$.
	\WHILE{$t < t^{max}$}
	\STATE - Compute $f_{cur} = f(\m{x})$
	\STATE  - Randomly generate search directions $\mathcal{D}$
	\STATE - Define $S = \argmin_{\m{d}\in\mathcal{D}}\{f(\m{x}+\delta \m{d}): \m{x}+\delta \m{d} \in \C{B}\}$
	\STATE - Randomly sample $\m{d}^*\in S$, define $f^*=f(\m{x}+\delta \m{d}^*)$
	\IF{$f^* < f^{\C{B}}_{min}$}
	\STATE $f^{\C{B}}_{min} \gets f^*, \m{x}_{min}^{\C{B}}\gets \m{x}+\delta \m{d}^*$
	\ENDIF
	\IF{$f^* > f_{cur}$}
	\STATE $\delta \gets \min\{\tau\delta,\delta_{max}\}$ 
	\ELSIF{$f^* < f_{cur}$} 
	\STATE  $\delta \gets \max\{\delta/\tau,\delta_{min}\}$
	\ENDIF
	\STATE$\m{x}\gets \m{x}+\delta \m{d}^*$
	\STATE $t \gets t + |\C{D}| + 1$
	\ENDWHILE
	\STATE \textbf{return}  $f^{\C{B}}_{min}, \m{x}^{\C{B}}_{min}$.
	\end{algorithmic}
\end{algorithm}
We provide two different options to generate the search directions $\mathcal{D}$: 
\begin{itemize}
    \item [(a)] By coordinate strategy, for each axis $i$, we take two directions $\m{e}_i$ and $-\m{e}_i$ with a probability $w_i$ calculated from variable importance.    
    \item [(b)] By random sampling from the unit sphere in $\R^p$.
\end{itemize}
The first option is more suitable when the discontinuous boundaries are in parallel to the axes, for example, ensemble tree regression functions. The second strategy works for general stepwise functions when the boundaries for each level set is not axis-parallel. We store the feasible sampled points $\m{x} + \delta \m{d}$ with their function values generated during local search in order to update the best function value $f_{\C{H}_i}$ for each hyper-rectangle $\C{H}_i$ for a new partition.

\subsection{Global Convergence of \texttt{StepDIRECT}}
Now, we are able to provide the global convergence result for \texttt{StepDIRECT}. We assume that $\m{w}=(1/p, \ldots, 1/p)\in \R^p$ for simplicity and the following theorem is still valid as long as $w_j>0$ for every $j\in \{1,\ldots,p\}$.

\begin{theorem}\label{thm:global_convergence_direct_v}
Suppose that $\m{w}=(1/p, \ldots, 1/p)\in \R^p$ and $f$ is continuous in a neighborhood of a global optimum. Then, \texttt{StepDIRECT} converges to the globally optimal function value for the stepwise function $f$ over the bounded box defined in~\eqref{eq:original_problem_1}.
\end{theorem}

\section{Numerical Experiments}\label{sec:numerical_examples}
In this section, we test the performance of \texttt{StepDIRECT} algorithm on two problems: optimization for random forest regression function and hyper-parameter tuning for classification. As explained in Section \ref{sec:motivating_examples}, we need to minimize a stepwise target function. We denote \texttt{StepDIRECT-0} by a variant of \texttt{StepDIRECT} when we skip Local\_Search in Step a of Algorithm \ref{alg:stepdirect}. 

\subsection{Optimization for Random Forest Regression Function}
We consider the minimization for Random Forest regression function over a bounded box constraint. We used the \textit{boston}, \textit{diabetes}, \textit{mpg} and  \textit{bodyfat}  data sets from the UCI Machine Learning Repository \cite{Dua:2019}. Table \ref{table:dataset_rf} provides the details of these four data sets.
\begin{table}[ht]
\scriptsize
\begin{center}
\caption{Dataset statistics, $p$: number of features, $N$: number of samples}\label{table:dataset_rf} 
\begin{tabular}{ |c|c|c|c|c|}
\hline
 Dataset  & boston & diabetes & mpg & bodyfat\\\hline
 $N$ & 506& 768 &234 & 252 \\ \hline
 $p$ & 13 & 8 & 11 & 14 \\ \hline
\end{tabular}
\end{center}
\end{table}

We train the random forest regression function on these data sets with 100 trees and use default settings for other parameters in \texttt{scikit-learn} package \cite{pedregosa2011scikit}. For comparison, we run the following optimization algorithms: \texttt{DIRECT} \cite{jones1993lipschitzian}\footnote{https://scipydirect.readthedocs.io/en/latest/}, \texttt{DE} \cite{price2006differential}\footnote{https://docs.scipy.org/doc/scipy/reference/index.html}, \texttt{PSO} \cite{Kennedy95}\footnote{https://pypi.org/project/pyswarms/}, \texttt{RBFOpt}   \cite{costa2018rbfopt}\footnote{https://projects.coin-or.org/RBFOpt}, \texttt{StepDIRECT-0}, and \texttt{StepDIRECT}.  For both \texttt{StepDIRECT-0} and \texttt{StepDIRECT} in Algorithm \ref{alg:stepdirect}, we set the maximum number of function evaluations $m^{max}=2000$. The same function evaluation limit is used for \texttt{DIRECT}. In Algorithm \ref{alg:LS}, we select $\delta = 1, \delta_{min}=0.001, \delta_{max}=2.5$, $\tau=1.5, |\C{D}| = 5$, and $t^{max}=1.5p$.
For all other algorithms, we use the default settings.
We run all algorithms for 20 times and report the mean and standard deviation results for final objective function values (denoted by ``obj") and running times in seconds (denoted by ``time") in Table \ref{table:opt_rf}.

\begin{table}[t!]
\scriptsize
\begin{center}
\caption{Optimization of random forest regression functions for four data sets (mean$\pm$standard deviation). ``obj" is the final objective function value, ``time" is the running time in seconds. An entry is in bold if the mean value is the lowest in the column. 
}\label{table:opt_rf} 
\begin{tabular}{|@{\hspace{0.06cm}}C{1.5cm}|l@{\hspace{0.08cm}}|c@{\hspace{0.04cm}}|@{\hspace{0.08cm}}c@{\hspace{0.04cm}}|@{\hspace{0.08cm}}c@{\hspace{0.04cm}}|@{\hspace{0.04cm}}c@{\hspace{0.04cm}}|}
\hline
          \multicolumn{2}{|c|}{Algorithms}      & boston      & diabetes     & mpg         & bodyfat     \\[1mm]
          \hline
\multirow{4}{*}{\texttt{DIRECT}} & obj  & 32.60 & 77.66  & 10.55  & 1.31   \\
                           &   &  ($\pm$ 0.00) &  ($\pm$ 0.00)  & ($\pm$ 0.00) & ($\pm$ 0.00)  \\
                           & time & 8.21  & 9.14   & 1.52  & 6.15   \\
                           &  & ($\pm$ 0.01)  & ($\pm$ 0.06)   & ($\pm$ 0.04)  & ($\pm$0.02)  \\ \hline
\multirow{4}{*}{\texttt{StepDIRECT-0}} & obj  & 28.66 & 75.8  & 10.37 & 1.27  \\
                          &   &  ($\pm$ 0.00) &  ($\pm$ 0.00)  &  ($\pm$ 0.00) & ($\pm$ 0.00)  \\
                           & time & 7.59   & 9.41   & 1.66  & 6.03   \\
                           &  &  ($\pm$ 0.01)  & ($\pm$ 0.16)   &  ($\pm$ 0.09)  & ($\pm$ 0.03)  \\ \hline
\multirow{4}{*}{\texttt{StepDIRECT}} & obj  & \textbf{28.35} & \textbf{58.68}  & 9.85 & \textbf{1.26}  \\ 
                           &   & \textbf{($\pm$0.02)} & \textbf{($\pm$ 1.04)}  & ($\pm$ 0.03)  & \textbf{($\pm$ 0.01)}  \\
                           & time & 7.87  & 10.74 & 1.88  & 6.63 \\
                           &  & ($\pm$ 0.23)  & ($\pm$ 0.15)  & ($\pm$ 0.17)  & ($\pm$ 0.26)  \\\hline
\multirow{4}{*}{\texttt{DE}}        & obj  & 28.40 & 69.71   & \textbf{9.82}  & 1.29  \\
                           &   & ($\pm$ 0.08) & ($\pm$ 1.72)  & \textbf{($\pm$ 0.01)}  & ($\pm$ 0.01)   \\
                           & time & 15.22  & 23.86  & 8.85  & 3.34   \\
                           &  & ($\pm$ 0.79)  & ($\pm$ 0.45)  & ($\pm$ 0.70)  & ($\pm$ 0.29)  \\ \hline
\multirow{4}{*}{\texttt{PSO}}     & obj  & 28.41 & 80.55 & 9.87 & 1.28  \\
                           &   & ($\pm$ 0.08) & ($\pm$ 3.82)  & ($\pm$ 0.06)  & ($\pm$ 0.05)  \\
                           & time & 30.60  & 34.67  & 34.47  & 29.49  \\
                           &  & ($\pm$ 0.49) &  ($\pm$ 11.91) &  ($\pm$ 0.88) &  ($\pm$ 0.14) \\\hline
\multirow{4}{*}{\texttt{RBFOpt}}    & obj  & 28.56  & 87.34  & 10.40  & 1.33   \\
                           &   &  ($\pm$ 0.02) &  ($\pm$ 3.82)  & ($\pm$ 0.00) &  ($\pm$ 0.04)  \\
                           & time & 14.70  & 11.95  & 8.70  & 12.83 \\
                           &  & ($\pm$ 2.85) & ($\pm$ 1.17)  &  ($\pm$ 0.08)  & ($\pm$ 1.83)\\\hline
\end{tabular}
\end{center}
\end{table}

From Table \ref{table:opt_rf}, we see that the function values returned by \texttt{StepDIRECT-0} are better than those of the original \texttt{DIRECT}. It illustrates the benefit of our proposed strategy for identifying potentially optimal hyper-rectangles and dividing hyper-rectangles which efficiently exploits the stepwise function structure. For \texttt{DIRECT} and \texttt{StepDIRECT-0}, their objective function outputs do not change for different runs since they are deterministic algorithms.         

To further see the advantage of local search incorporated into \texttt{StepDIRECT-0}, in \texttt{StepDIRECT} the local search is initialized with the best point in the potentially optimal hyper-rectangle and runs with search directions $\mathcal{D}$ randomly generated by the coordinate strategy. Compared with \texttt{StepDIRECT-0}, we notice that the \texttt{StepDIRECT} algorithm achieved lower objective function values. From Table \ref{table:opt_rf}, \texttt{StepDIRECT} shows the best overall performance in terms of solution quality except for \textit{mpg} dataset. By embedding the local search, we can significantly improve the solution quality. In general, the proposed algorithm runs faster than other baseline methods \texttt{DE}, \texttt{PSO}, and \texttt{RBFOpt}, except for the run time of \textit{bodyfat} dataset, \texttt{DE} outperforms \texttt{StepDIRECT}.

\subsection{Hyper-parameter Tuning for Classification}
We tune the hyper-parameters for: 1) multi-class classification with linear support vector  machines (SVM) and logistic regression, and 2) imbalanced binary classification with RBF-kernel SVM \cite{hastie2009elements}. 
We use three datasets: \textit{MNIST} from \cite{MNIST}, \textit{PenDigits} from the  UCI  Machine  Learning  Repository \cite{Dua:2019}, and a synthetic dataset \textit{Synth} generated by the Mldatagen generator \cite{TOMAS2014155}. For all datasets, we set the ratio among training, validation and test data partitions as $3:1:1$ and report the performance on the test dataset.  



\subsubsection{Multi-class Classification}

For multi-class classification problems with the number of classes $K \ge 3$, there are two widely-used approaches:  “one-vs-all” and “one-vs-one” strategies.

\textbf{One-vs-all classification:} In this experiment, we tune the hyper-parameters for multi-class support vector machines and logistic regression. 
For multi-class classification problems, we take the ``one-vs-all" approach to turning a binary classifier into a multi-class classifier. 
 For each class $k\in \{1,\ldots,K\}$ and a given hyper-parameter $C_k$, the one-vs-all SVM solves the following problem for the dataset $\{\m{x}_i,y_i\}_{i=1}^N$
\begin{align}\label{eq:one-vs-all SVM}
 \min_{\{\m{w}^k, b^k\}} & \dfrac{1}{2}\norm{\m{w}^k}^2 + C_k \sum_{i=1}^{N}\xi_{i}^k & \\\nonumber
 \text{s.t. } & (\m{w}^k)^T \m{x}_i + b^k \geq 1 - \xi_{i}^k, & \quad \text{if } y_i = k\\\nonumber
 & (\m{w}^k)^T\m{x}_i + b^k \leq -1 + \xi_{i}^k, &\quad \text{if } y_i \neq k\\\nonumber
 & \xi_{i}^k \geq 0,\, i\in \{1,\ldots,N\}. \nonumber
\end{align}
The class of each point $\m{x}_i$ is determined by 
\begin{equation}
\text{class of } \m{x}_i =  \argmax_{k\in \{1,\ldots,K\}} \{(\m{w}^k)^T\m{x}_i + b^k \}.
\end{equation}
Different from many default implementations by taking the same value for all $C_k, 1\le k \le K$, we allow them to take different values because the margin between different classes may vary significantly from each other. The one-vs-all logistic regression follows the same approach by replacing \eqref{eq:one-vs-all SVM} with the binary logistic regression classifier.

We search $K$ hyper-parameters $C_k$ for $K$ classifiers in the log space of $\Omega = [10^{-3}, 10^3]^{K}$. The number of hyper-parameters for each dataset is given in Table \ref{tab:one_vs_all_svm}, which is denoted by ``$p$".  For comparison, we run the following algorithms: \texttt{Random Search (RS)}, \texttt{DIRECT}, and \texttt{Bayesian Optimization (BO)\footnote{https://rise.cs.berkeley.edu/projects/tune/}}. The widely used \texttt{Grid Search (GS)} is not considered here because \texttt{GS} is generally not applicable when the search space dimension is beyond 5. For all algorithms, the computational budget is $m^{max}= 100K$ in terms of the number of training $K$ base classifiers. The results for one-vs-all SVM and one-vs-all logistic regression are shown in Tables \ref{tab:one_vs_all_svm} and \ref{tab:one_vs_all_lr}, respectively. 
We can observe that \texttt{StepDIRECT} gives the best test accuracy for these data sets. Compared with the random search \texttt{RS}, \texttt{StepDIRECT} improves the test accuracy $0.6-2.0\%$. We notice that the original \texttt{DIRECT} algorithm makes little improvement and often gets stuck in the local region until consuming all running budgets, while the \texttt{StepDIRECT} algorithm can make consistent progress by balancing the local exploration and global search. 

\begin{table}[ht]
\scriptsize
\begin{center}
\caption{One-vs-all SVM, $K$: the number of classes, $p$: the number of tuning parameters}\label{tab:one_vs_all_svm}
\begin{tabular}{ |c|c|c|c|c|c|c| } 
 \hline
 Dataset & $K$ & $p$ & \texttt{RS} & \texttt{DIRECT} & \texttt{StepDIRECT} & \texttt{BO}\\\hline
 Synth & 3 & 3 & 74.8\% & 75.3\% & \textbf{76.8}\% & 75.0\%\\ \hline
 PenDigits & 10 & 10 & 94.0\% & 94.2\%  & \textbf{95.2}\% & 94.9\% \\\hline
 MNIST & 10 & 10  & 91.3\% & 91.8\% & \textbf{92.2}\% & 92.0\% \\\hline
\end{tabular}
\end{center}
\end{table}

\begin{table}[ht]
\scriptsize
\begin{center}
\caption{One-vs-all logistic regression, $K$: the number of classes, $p$: the number of tuning parameters}\label{tab:one_vs_all_lr}
\begin{tabular}{ |c|c|c|c|c|c|c| } 
 \hline
 Dataset & $K$ & $p$ & \texttt{RS} & \texttt{DIRECT} & \texttt{StepDIRECT} & \texttt{BO}\\\hline
Synth & 3 & 3 & 74.5\% & 75.2\% & \textbf{75.5}\% & \textbf{75.5}\%\\ \hline
 PenDigits & 10 & 10 & 94.8\% & 95.0\% & \textbf{95.3}\% & 95.1\%\\\hline
 MNIST & 10 & 10 & 92.0\% & 92.2\% & \textbf{92.6}\% & 92.3\%\\\hline
\end{tabular}
\end{center}
\end{table}

\textbf{One-vs-one classification:} For each pair of classes $i$ and $j, 1\leq i < j \leq K$, we need to tune an associated model parameter $C_{ij}$. The one-vs-one SVM solves the following problem
\begin{align}\label{eq:one-vs-one SVM}
 \min_{\{\m{w}^{ij}, b^{ij}\}} & \dfrac{1}{2}\norm{\m{w}^{ij}}^2 + C_{ij} \sum_{t=1}^{N}\xi_{t}^{ij} & \\\nonumber
 \text{s.t. } & (\m{w}^{ij})^T\m{x}_t + b^{ij} \geq 1 - \xi_{t}^{ij}, & \quad \text{if } y_t = i\\\nonumber
 & (\m{w}^{ij})^T\m{x}_t + b^{ij} \leq -1 + \xi_{t}^{ij}, &\quad \text{if } y_t = j\\\nonumber
 & \xi_{t}^{ij} \geq 0, \, t\in \{1,\ldots,N\} \nonumber
\end{align}
for all pairs $(i,j), 1\leq i < j \leq K$. There are $\frac{K(K-1)}{2}$ hyper-parameters $C_{ij}$ to be tuned.
We use the following voting strategy: if $\sign((\m{w}^{ij})^T\m{x}_t + b^{ij})$ says $\m{x}_t$ is in the $i$-th class, then the vote for the $i$-th class is added by one. Otherwise, the vote for the $j$-th class is increased by one. Then, the class of $\m{x}_t$ has the largest vote. Similar to the one-vs-all case, we search different hyper-parameters for all pair classifiers ($\Omega = [10^{-3}, 10^3]^{K(K-1)/2}$). For all algorithms, the budget is $10K(K-1)$ in terms of training a base classifier. The results are shown in Tables \ref{tab:one_vs_one_svm} and \ref{tab:one_vs_one_lr} for one-vs-one SVM and one-vs-one logistic regression, respectively. Overall, \texttt{StepDIRECT} performs the best in many cases compared to the base-line algorithms. Especially, it outperforms the Bayesian optimization algorithm \texttt{BO} in 4 out of 6 cases. 

\begin{table}[ht]
\scriptsize
\begin{center}
\caption{One-vs-one SVM, $K$: the number of classes, $p$: the number of tuning parameters}\label{tab:one_vs_one_svm}
\begin{tabular}{ |c|c|c|c|c|c|c| } 
 \hline
 Dataset & $K$ & $p$ & \texttt{RS} & \texttt{DIRECT} & \texttt{StepDIRECT} & \texttt{BO}\\\hline
 Synth & 3 & 3 & 78.0\% & 78.5\% & \textbf{79.5}\% & 78.8\%\\ \hline
 PenDigits & 10 & 45 & 97.9\% & 98.2\% & \textbf{98.8}\% & 98.6\%\\\hline
 MNIST & 10 & 45 & 91.5\% & 92.2\% & 93.0\%  & \textbf{93.2}\%\\\hline
\end{tabular}
\end{center}
\end{table}

\begin{table}[ht]
\scriptsize
\begin{center}
\caption{One-vs-one logistic regression, $K$: the number of classes, $p$: the number of tuning parameters}\label{tab:one_vs_one_lr}
\begin{tabular}{ |c|c|c|c|c|c|c| } 
 \hline
 Dataset & $K$ & $p$ & \texttt{RS} & \texttt{DIRECT} & \texttt{StepDIRECT} & \texttt{BO}\\\hline
 Synth & 3 & 3 & 77.5\% & 77.75\% & 78.25\%& \textbf{78.5}\%\\ \hline
 PenDigits & 10 & 45 & 96.8\% & 97.8\% & \textbf{98.4}\%  & 97.4\% \\\hline
 MNIST & 10 & 45 & 93.3\% & 93.8\% & \textbf{94.1}\%  & 93.7\% \\\hline
\end{tabular}
\end{center}
\end{table}

\subsubsection{Imbalanced Binary Classification}
In the second example, we consider to tune parameters of RBF-SVM with $\ell_2$ regularization $C$, kernel width $\gamma$, and class weight $c_p$ for the imbalanced binary classification problems. In this setting, we compare the performance of \texttt{StepDIRECT} with \texttt{DIRECT}, \texttt{RS}, \texttt{BO}, and \texttt{GS} in tuning a three-dimensional hyper-parameter $w=(C, \gamma, c_p)$ to achieve a high test accuracy. 

In our experiments, we used five binary classification datasets, as shown in Table \ref{tab:imbalanced_datasets}, which can be obtained from LIBSVM package \cite{CC01a}.
\begin{table}[ht]
\scriptsize
\begin{center}
\caption{Data set statistics, $p$: the number of features, $N$: the number of samples, $N_{-}/N_{+}$: the class distribution ratio}\label{tab:imbalanced_datasets}
\begin{tabular}{ |c|c|c|c| } 
 \hline
 Dataset  &   $p$ & $N$ &  $N_{-}/N{+}$\\\hline
 fourclass & 2 & 862 & 1.8078 \\\hline
diabetes & 8 & 768 & 1.8657  \\\hline
statlog & 24 & 1,000 & 2.3333  \\\hline
svmguide3 & 22 & 1,284 & 3.1993\\\hline
ijcnn1 & 22 & 141,691 & 9.2643 \\\hline
\end{tabular}
\end{center}
\end{table}
\begin{table}[t]
\scriptsize
\begin{center}
\caption{RBF-SVM} 
\label{tab:RBF_SVM}
\begin{tabular}{ |c@{\hspace{1mm}}|c@{\hspace{1mm}}|c@{\hspace{1mm}}|c@{\hspace{1mm}}|c@{\hspace{1mm}}|c| } 
 \hline
 Dataset  & \texttt{RS} &  \texttt{GS} & \texttt{DIRECT} & \texttt{StepDIRECT} & \texttt{BO}\\\hline
{\scriptsize diabetes} & 79.8\% & 81.2\% & 79.2\% & \textbf{81.8}\% &81.3\%  \\\hline
 {\scriptsize fourclass}  & \textbf{100.0}\% & 80.9\% & \textbf{100.0}\% & \textbf{100.0}\% & \textbf{100.0}\%\\\hline
 {\scriptsize statlog}  & 77.5\% & 78.5\% & 80.5\% & \textbf{81.6}\% & 81.0\%\\\hline
 {\scriptsize svmguide3}  & 83.5\% & 83.9\% & \textbf{84.3}\% & \textbf{84.3}\% & 82.9\%\\\hline
  {\scriptsize ijcnn1}  & 97.6\% & \textbf{98.3}\% & 94.6\% & 98.2\% & \textbf{98.3}\%\\\hline
\end{tabular}
\end{center}
\end{table}

For all algorithms, the feasible set is chosen as $C\in [10^{-3}, 10^{3}]$, $\gamma \in [10^{-6}, 10^0]$ and $c_p\in [10^{-2}, 10^2]$. The search is conducted in the log-scale. For the \texttt{GS} algorithm, we uniformly select 5 candidates for each parameter. For a fair comparison, we set the budget as 125 for the limit of number of function evaluations for all algorithms. Table \ref{tab:RBF_SVM} shows the test accuracy for the experiment. 



As we can see from Table \ref{tab:RBF_SVM}, \texttt{StepDIRECT} is always the top performer when compared with  \texttt{DIRECT} and \texttt{RS} algorithms. Furthermore, the proposed algorithm achieves competitive performance with \texttt{BO}.

\section{Conclusion}\label{sec:conclusion}
In this paper, we have proposed the \texttt{StepDIRECT} algorithm for solving the black-box stepwise function, which can exploit the special structure of the problem. We have introduced a new definition to identify the  potentially optimal hyper-rectangles and divide the hyper-rectangles. A stochastic local search is embedded to improve solution quality and speed up convergence. The global convergence is proved and numerical results on two practical machine learning problems show the state-the-art-performance of algorithm. In the future, we plan to combine the ideas from \texttt{StepDIRECT} and Bayesian Optimization, and develop practical algorithms for tuning  discrete hyper-parameters in machine learning and deep neural networks.


\bibliography{reference}
\bibliographystyle{plain}


\newpage
\onecolumn
\appendix
\section{Supplementary Material}

\textbf{Lemma \ref{lem_0001}}. \textit{Let $\epsilon>0$ be the positive constant used in Definition \ref{def:new_potential_optimal} and $f_{min}$ be the current best function value. Let $\mathcal{I}$ be the set of indices of all existing hyper-rectangles. For each $j\in \mathcal{I}$, let $
\mathcal{I}_l = \{i\in \mathcal{I}: d_i \sigma_i < d_j \sigma_j\}, 
\mathcal{I}_b = \{i\in \mathcal{I}: d_i \sigma_i > d_j \sigma_j\},
\mathcal{I}_e = \{i\in \mathcal{I}: d_i \sigma_i = d_j \sigma_j\}
$
and $g_{i} = \frac{f_{\C{H}_i} - f_{\C{H}_j}}{d_i \sigma_i - d_j \sigma_j}$ for all $i\neq j$.
The hyper-rectangle $\C{H}_j$ is potentially optimal if and only if three following conditions hold:
\begin{itemize}
    \item[(a)] 
$f_{\C{H}_j} \leq f_{\C{H}_i}$ for every $i \in \mathcal{I}_e;$
\item[(b)] 
$
\max_{i\in \mathcal{I}_l}g_{i}\leq  \min_{i\in \mathcal{I}_b}g_{i};
$
\item[(c)] If $f_{median} > f_{min}$ then 
\begin{align*}
    \epsilon \leq \frac{f_{min} - f_{\C{H}_j}}{|f_{min}-f_{median}|}  + \frac{d_j \sigma_j \min_{i\in \mathcal{I}_b}g_{i}}{|f_{min}-f_{median}|}; 
\end{align*}
otherwise,  
\begin{equation*}
f_{\C{H}_j} \leq d_j \sigma_j \min_{i\in \mathcal{I}_b}g_{i} + f_{min}.
\end{equation*}
\end{itemize}}

\begin{proof}
First, assume that $\C{H}_j$ is potentially optimal. For $i\in\mathcal{I}_e$, the inequality $f_{\C{H}_j} \leq f_{\C{H}_i}$ follows directly from \eqref{eq:poc_1}. 

For $i\in \mathcal{I}_l$, from  \eqref{eq:poc_1},  we have
\[
K\geq \frac{f_{\C{H}_j} - f_{\C{H}_i}}{d_j\sigma_j - d_i\sigma_i},
\]
and for $i\in \mathcal{I}_b$, it implies that
\[
K \leq \frac{f_{\C{H}_i} - f_{\C{H}_j}}{d_i\sigma_i - d_j \sigma_j}.
\]
Hence, (b) directly follows from above by taking the maximum over $\mathcal{I}_l$ and taking the minimum over $\mathcal{I}_b$.

By \eqref{eq:poc_2}, when $f_{min}\neq f_{median}$, we have
\[
\epsilon\leq \frac{f_{min} - f_{\C{H}_j}}{|f_{min}-f_{median}|} + K\frac{d_j \sigma_j}{|f_{min}-f_{median}|}.
\]
Eq. \eqref{eq:lemma1_3} is a consequence of the above inequality by taking
\[K=\min_{i\in \mathcal{I}_b}\frac{f_{\C{H}_i} - f_{\C{H}_j}}{d_i \sigma_i - d_j \sigma_j}.\]
Similar arguments hold for \eqref{eq:lemma1_4} when $f_{min}=f_{median}$.

For the other direction of the proof, from conditions (a) and (b), we can derive \eqref{eq:poc_1}. Using conditions (b) and (c) we can get \eqref{eq:poc_2}.  
\end{proof} 


\textbf{Theorem \ref{thm:global_convergence_direct_v}}. \textit{Suppose that $\m{w}=(1/p, \ldots, 1/p)\in \R^p$ and $f$ is continuous in a neighborhood of a global optimum. Then, \texttt{StepDIRECT} converges to the globally optimal function value for the stepwise function $f$ over the bounded box defined in~\eqref{eq:original_problem_1}.}

\begin{proof}
We will prove that for any distance $\delta > 0$ to the optimal solution $\m{x}^*$, \texttt{StepDIRECT} will sample at least a point within a distance $\delta$ of $\m{x}^*$. We will argue that the points sampled by the algorithm  form a dense subset of the unit hypercube $\Omega$.

We note that the new hyper-rectangles are generated by splitting exiting ones into one thirds on some dimensions.   
The side lengths for one hyper-rectangle are in of the form $3^{-k}$ for some $k \ge 0$. Notice that as we always divide the larger sides, we have to divide those sides with length $3^{-k}$  before dividing any side of length $3^{-(k+1)}$. After carrying out $r$ divisions, the hyper-rectangle will have $j=\mod(r,p)$ sides of length $3^{-(k+1)}$ and $n-j$ sides of length $3^{-k}$ with $k=(r-j)/p$. Hence, the diameter (the largest distance from one vertex to the other vertex) of the hyper-rectangle is given by 
\begin{equation}\label{eq: diameter_of_hyper_rectangle}
d = [j3^{-2(k+1)} + (n-j)3^{-2k}]^{0.5},
\end{equation}
and the volume is $V=3^{-(nk+j)} = 3^{-r}$.
We can see that, as the number of divisions approaches to $\infty$, both the diameter $d$ and the volume $V$ converge to 0.

At the $t$-th iteration of \texttt{StepDIRECT}, let $l_t$ be the fewest number of partitions undergone by any hyper-rectangle.
It follows that these hyper-rectangles would have the largest diameter. We will show by contradiction that  
\[
\lim_{t\rightarrow\infty} l_t =\infty.
\]
If $\lim_{t\rightarrow\infty} l_t$ is bounded with an upper limitation $r < \infty$, we define the following nonempty set 
\[
\mathcal{I}_{r} = \{i: V_{i} = 3^{-r}\}.
\]
We choose a hyper-rectangle $\C{H}_j$ where $j\in J:= \argmax_{i}\{d_i \sigma_i\}$ with the best function value $f_{\C{H}_j}$ over the hyper-rectangle $\C{H}_j$. By the update rule, this hyper-rectangle will be potentially optimal as the conditions in Definition \ref{def:new_potential_optimal} are fulfilled with $\tilde{K}>\max\{K_1, K_2\}$, where 
\begin{align}
&K_1 = \dfrac{f_{\C{H}_j}- f_{min} + \epsilon|f_{min}-f_{median}|}{d_j \sigma_j}\\\nonumber
&K_2 = \max \{\dfrac{f_{\C{H}_j} - f_{\C{H}_i}}{d_j \sigma_j - d_i \sigma_i}: d_j \sigma_j \neq d_i \sigma_i\}. \nonumber
\end{align}
If $j\in\mathcal{I}_{r}$, it follows that $|\mathcal{I}_{r}|$ decreases by 1. Otherwise, a series of new hyper-rectangles will be generated with volumes no greater than $\frac{V_j}{3}$. We can repeat the above process for infinite times. Notice that the variability $\sigma_j$ is bounded within $[\epsilon_\sigma, 1]$. After a finite time of iterations, at least one hyper-rectangle $j\in \mathcal{I}_{r}$ will be selected. It follows $\mathcal{I}_{r}$ would eventually diminish to an empty set. This contradiction proves that $\lim_{t\rightarrow\infty} l_t =\infty$. As a result, the points generated by \texttt{StepDIRECT} will be dense in the original rectangle $\Omega$ and the global convergence directly follows from the continuity assumption.
\end{proof}

\end{document}